\documentclass[10pt]{article}
\UseRawInputEncoding
\usepackage{amssymb}
\usepackage{amsfonts}
\title{Countable homogeneous ordered bipartite graphs}
\author{J K Truss} 
\date{University of Leeds}

\begin{document}
\maketitle  
\newtheorem{lemma}{Lemma}[section]
\newtheorem{theorem}[lemma]{Theorem}
\newtheorem{corollary}[lemma]{Corollary}

\setcounter{footnote}{1}\footnotetext{
keywords: bipartite graph, homogeneous, ordered}
\newcounter{number}

\begin{abstract} A classification is given of all the countable homogeneous ordered bipartite graphs
\end{abstract}

\section{Introduction}

A relational structure is said to be homogeneous (also called `ultrahomogeneous', as in \cite{Lachlan1} for instance) 
if any isomorphism between finite substructures extends to an automorphism of the whole structure. The study of
homogeneous structures as highly symmetrical and beautiful mathematical constructs in their own right has centred
on classifications which exist in now rather a large number of individual cases. Relations with other areas were
however explored by Ne$\breve{\rm s}$et$\breve{\rm r}$il and R\"odl in \cite{Nesetril}, where interesting connections 
with Ramsey theory were found. A restriction for their work is that the structures should carry a linear 
order. That is partly the justification for the present piece of work, in which we present a classification of 
the countable homogeneous ordered bipartite graphs.

The classification of countable (unordered) homogeneous bipartite graphs is rather straightforward, and is recalled for instance 
in \cite{Goldstern, Jenkinson}. There are just five kinds, empty, complete, perfect matching and its complement, and generic. The 
first four exist in all finite and countable cardinalities, the last just in the infinite case, and is characterized by 
the property that for any two finite disjoint subsets in one of the parts, there is member of the other part (a `witness') 
joined to all elements of the first set and to none of the second. We at once remark that to get an analogous notion in 
the ordered case, in which the two parts will be ordered like the rationals, one has to say that the set of witnesses 
must be dense, and the principal difficulties in this piece of work are in ensuring that this holds.

There is a slight ambiguity in the notion of bipartite graph, as to whether the two parts should be distinguished (by unary
predicates). This amounts to knowing whether or not the automorphism group can interchange the two parts or not. In the unordered
cases one can always do this (whereas in \cite{Goldstern} this is not possible, since the two parts have unequal cardinalities).
For ease therefore, we shall throughout designate our two parts as `red' and `blue'.

What we do not do here, and it is not clear how difficult, or fruitful it might be, is to extend these results to the multipartite case,
as was done in the unordered case in \cite{Jenkinson}. We leave this for later study, but make some brief remarks on it at the end.

I would like to thank David Bradley-Williams for introducing me to this problem.

\vspace{.1in}

\noindent{\bf Outline of the paper}

As just explained, we shall start with a bipartite graph $(X, R)$ which also carries a linear ordering $<$, and the two parts will 
be indicated by colours red and blue. The first stage is to list the possibilities for the reducts of $(X, <, R)$, viewed just as
coloured linear orders. Then in the main result, these possibilities are examined in turn to see what (homogeneous) graph structures
are compatible with each. There are some `degenerate' cases, and three cases of some genuine interest in which the coloured linear order
is two copies of the rationals $(2.{\mathbb Q}, <)$, a rational family of two-element sets $({\mathbb Q}.2, <)$, or where it is the 
`2-coloured version of the rationals $({\mathbb Q}_2, <)$', which is a countable dense linear order without endpoints each point
coloured `red' or `blue', so that the colours both occur densely.

The most involved cases are the `generic' ones, in which one has to verify that witnesses as mentioned earlier arise densely. 

\section{Subdivision into cases}

\begin{lemma} Let $(X, <, R)$ be a countable homogeneous ordered bipartite graph in which $(X, <)$ is a linear order, and
$(X,R)$ is a bipartite graph with parts distinguished by `red' and `blue'. Then as a 2-coloured linear ordering, $(X, <)$ is  isomorphic to one of the following:

\noindent (I) Degenerate cases

(i) a singleton, or a doubleton each point of colour red or blue (and one of each),

(ii) $({\mathbb Q}, <)$ monochromatically coloured (red or blue),

(iii) $({\mathbb Q} \cup \{\infty\}, <)$ or $({\mathbb Q} \cup \{-\infty\}, <)$ with $\mathbb Q$ coloured by one of the colours,
and the singleton coloured by the other,,

\noindent (II) Non-degenerate cases

(iv) $(2.{\mathbb Q}, <)$ with one copy of $\mathbb Q$ coloured by one of the colours, the other by the other,

(v) $({\mathbb Q}.2, <)$ with each (adjacent) pair coloured the same, either red-blue, or blue-red,

(vi) the generic 2-coloured linear order ${\mathbb Q}_2$.
 \label{2.1}  \end{lemma}

\noindent{\bf Proof}: Let $G = {\rm Aut}(X, <, R)$. Since on any monochromatic subset, the bipartite graph relation 
is empty, $G$ acts transitively on each of the sets of red and blue singletons, and also on the sets of red and blue doubletons. Let $X_r$, $X_b$
be the sets of red, blue points. It follows that each of $X_r$, $X_b$ is empty, or is a singleton, or is isomorphic to 
$\mathbb Q$. 

If $X$ is finite, we get cases (i).

If $X$ is monochromatic non-singleton, we get cases (ii).

If $X$ is infinite, but one of $X_r$, $X_b$ is a singleton, we get cases (iii), noting that the singleton $z$ cannot lie between 
points $x$, $y$ of the other colour, since an automorphism taking $x$ to $y$ would have to move $z$.

From now on, we assume that $X_r$ and $X_b$ are both infinite, hence isomorphic to $\mathbb Q$. If $X_r < X_b$ 
or $X_b < X_r$ we get cases (iv). Otherwise there are $x < y < z$ where $x$ and $z$ have the same colour and $y$ the other colour.
Taking $x$ to $z$ by an automorphism and considering the image of $y$ we see that there are two red points with a blue
point in between, and also two blue points with a red point in between.

If some red point has an immediate successor, then it must be blue, and by transitivity of $G$ on $X_r$, every red point 
has an immediate successor, so we are in case (v) with each consecutive pair coloured (red, blue). Similarly if
some blue point has an immediate successor.

Now assume that none of the above hold. To show that $(X, <) \cong ({\mathbb Q}_2, <)$ it suffices to show that if $x < y$ 
then there are red and blue points in between. We already know this if $x$ and $y$ have the same colour. Now suppose that $x$ is red and $y$ is blue. 
Since $y$ is not the immediate successor of $x$, there is $z$ with $x < z < y$. If $z$ is red, there is a blue point 
between $x$ and $z$, and if $z$ is blue there is a red point between $z$ and $y$. $\Box$

\vspace{.1in}

\section{The main result}

Before we present the main classification, we remark on some of the structures which arise, which are of a generic flavour.
In view of the preliminary result given in the previous section, Lemma \ref{2.1}, we already have an outline of how we expect
the classification to go. The interesting cases correspond to the non-degenerate cases (iv) and (vi) in that list, and we first
give the constructions for these.

For case (iv) we start with the coloured linear order obtained by concatenating a `red' copy of the rationals ${\mathbb Q}_{red}$ and 
a blue copy ${\mathbb Q}_{blue}$, and we have to decide which red points to join to which blue points. One case of this, which we
refer to as `bounded generic' is obtained by starting with ${\mathbb Q}_2$, the $\{red, blue\}$ 2-coloured version of $\mathbb Q$, and
starting by defining a different ordering $\prec$ on it, given by $x \prec y$ if $x$ and $y$ have the same colour and $x < y$,
or $x$ is red and $y$ is blue. Thus we `separate out' the red and blue points, retaining their orders, and placing all red
points below all blue points. This is then clearly the same as ${\mathbb Q}_{red}$ followed by ${\mathbb Q}_{blue}$. The reason
for using the identification with ${\mathbb Q}_2$ is to use it to choose the bipartite graph relation $R$. Specifically, we
let $R(x,y)$ if $x$ is red and $y$ is blue and $x < y$, or $y$ is red and $x$ is blue and $y < x$. We call this `bounded',
since the set of neighbours of red $x$ is $(x, \infty) \cap {\mathbb Q}_{blue}$ (which is thus bounded below). There are
three other versions of this case, obtained by reversing the ordering of one or other colour, or both. In the {\em red reversal} 
${\mathbb Q}_{red}$ has its ordering reversed, but not ${\mathbb Q}_{blue}$, in its {\em blue reversal}, ${\mathbb Q}_{blue}$
has its ordering reversed, but not ${\mathbb Q}_{red}$, and in the red-blue reversal, both are reversed. In all these cases $R$
is unchanged. 

The other main case here is `unbounded generic', which is most easily obtained as the Fra\"iss\'e limit of the class
of all finite $\{red, blue\}$-coloured ordered bipartite graphs in which all red points precede all blue points. This will have underlying ordered 
set ${\mathbb Q}_{red}$ followed by ${\mathbb Q}_{blue}$, and for any finite disjoint subsets $A_1$ and $A_2$ of ${\mathbb Q}_{red}$,
the set of neighbours of all members of $A_1$ which are also non-neighbours of all members of $A_2$ is a dense subset of
${\mathbb Q}_{blue}$, and similarly for $B_1 \cup B_2 \subseteq {\mathbb Q}_{blue}$. A standard back-and-forth proof shows that
these statements are sufficient to characterize this case up to isomorphism.

For case (vi), the underlying ordered set is ${\mathbb Q}_2$, and we may view ${\mathbb Q}_{red}$ and ${\mathbb Q}_{blue}$ as subsets
(each dense in the whole order, unlike in case (iv)). Apart from the empty and complete cases, we have the `left complete' case
in which $R = \{(a, b): a \mbox{ red }, b \mbox{ blue }, a < b\}$ and `right complete',
$R = \{(a, b): a \mbox{ red }, b \mbox{ blue }, a > b\}$. In addition there are right and left generics, which are the Fra\"iss\'e limits
of the set of finite bipartite graphs on red/blue in which if $(a, b) \in R$ where $a$ is red and $b$ is blue, $a < b$,
and of those in which instead, $a > b$. There are also the fully generic, which is the union of these two, and also the union of the
left complete and right generic, and also the left generic and right complete. 

The right generic can in this case be characterized by saying that for any finite disjoint sets $A_1$, $A_2$ of red points, the set of
blue points which are related to all members of $A_1$ and to none of $A_2$, is dense in $(\max (A_1 \cup A_2), \infty)$, and for any
finite disjoint sets $B_1$, $B_2$ of blue points, the set of red points which are related to all members of $B_1$ and to none of $B_2$, 
is dense in $(- \infty, \min (B_1 \cup B_2))$. Again these conditions suffice for a uniqueness argument using back-and-forth. The
left generic can be similarly characterized. 

\begin{theorem} $(X, <, R)$ be a countable homogeneous ordered bipartite graph in which $(X, <)$ is a linear order, and
$(X,R)$ is a bipartite graph with parts distinguished by `red' and `blue'. Then $(X, <, R)$ is isomorphic to
one of the following, corresponding to the cases enumerated in Lemma \ref{2.1}:

\noindent (I) Degenerate cases

(i) a singleton, or a doubleton each point of colour red or blue (and one of each), with graph relation empty or complete,

(ii) $X = {\mathbb Q}$ under its usual ordering, monochromatically coloured (red or blue), and $R = \emptyset$,

(iii) $X = {\mathbb Q} \cup \{\infty\}$ or ${\mathbb Q} \cup \{-\infty\}$ with the usual ordering, where $\mathbb Q$ is coloured
by one of the colours, the singleton is coloured by the other, and $R$ is empty or complete.

\noindent (II) Non-degenerate cases

(iv) $(2.{\mathbb Q}, <)$ with one copy of $\mathbb Q$ coloured by one of the colours, the other by the other. Writing
these two copies as ${\mathbb Q}_{red}$ and ${\mathbb Q}_{blue}$, where we assume that ${\mathbb Q}_{red} < {\mathbb Q}_{blue}$, then
up to isomorphism the bipartite relation $R$ is empty or complete, or is bounded generic, possibly with one or both colours reversed,  
or unbounded generic,

(v) $({\mathbb Q}.2, <)$ with each (adjacent) pair coloured the same, either red-blue, or blue-red, in the former case
with these possibilities: $\{(a, b): b \mbox{ paired with } a\}$ (perfect matching),

$\{(a, b): b > \; \mbox{the blue point paired with } a\}$,

$\{(a, b): b < \; \mbox{the blue point paired with } a\}$,

or unions of these.

(vi) the generic 2-coloured linear order ${\mathbb Q}_2$, and the graph relation which is left complete, or 
right complete, or left generic or right generic, or unions of these.  \label{3.1}    \end{theorem}

\noindent{\bf Proof}: We remark that in all cases (i)-(vi), $R$ could be empty or complete, automatically giving rise to a homogeneous
structure, so from now on we disregard this. In particular, if $X$ is monochromatic, then $R = \emptyset$, and if $|X| = 2$, $R$ is
empty or complete. This deals with the cases (i) and (ii). Furthermore, in case (iii), if $R \neq \emptyset$, then by transitivity it
follows at once that it is complete. So we just have to deal with the non-degenerate cases (iv)-(vi).

(iv) Assume without loss of generality that ${\mathbb Q}_{red} < {\mathbb Q}_{blue}$. First we show that for each $a \in {\mathbb Q}_{red}$,
$|R(a)|, |{\mathbb Q}_{blue} \setminus R(a)| \neq 1$ (and similarly for $b \in {\mathbb Q}_{blue}$), where $R(a)$ is the set of neighbours
of $a$. By passing to the complement if necessary, it suffices to rule out the case in which $|R(a)| = 1$. If this holds, then by
1-transitivity the same holds for all $a \in {\mathbb Q}_{red}$. We show that $R$ must then be a perfect matching. Note that since there is
$b \in {\mathbb Q}_{blue}$ such that $R(a) = \{b\}$ for some $a \in {\mathbb Q}_{red}$, by transitivity, every $b \in {\mathbb Q}_{blue}$
has this form. If $R$ is {\em not} a perfect matching, then there must therefore be $a_1 < a_2$ in ${\mathbb Q}_{red}$ and 
$b \in {\mathbb Q}_{blue}$ such that $\{b\} = R(a_1) = R(a_2)$. Choose $b_3 \in {\mathbb Q}_{blue}$ not equal to $b$. By the above 
remark there is $a_3 \in {\mathbb Q}_{red}$ such that $R(a_3) = \{b_3\}$. Then $a_3 \neq a_1, a_2$ and so $a_3 < a_1 < a_2$ or
$a_1 < a_3 < a_2$ or $a_1 < a_2 < a_3$. In the first two cases let $p$ be given by $p(a_1) = a_3$, $p(a_2) = a_2$, and in the third
case let $p(a_1) = a_1$ and $p(a_2) = a_3$. Then $p$ is a finite partial automorphism, so by homogeneity extends to an automorphism
$\theta$. Since $p$ fixes $a_1$ or $a_2$, $\theta(b) = b$, and as $\theta(a_1)$ or $\theta(a_2)$ equals $a_3$, $\theta(b) = b_3$,
which gives a contradiction. We deduce that $R$ is a perfect matching.

Now choose $a_1 < a_2$ in ${\mathbb Q}_{red}$ and let $b \in {\mathbb Q}_{blue}$ be such that $R(a) = \{b\}$ for some $a$ between $a_1$ and $a_2$.
Let $p(a_1) = a_2$ and $p(b) = b$. This is a finite partial automorphism, but it cannot extend to any automorphism $\theta$, since
then $\theta(a)$ would have to equal $a$, contrary to $\theta(a_1) = a_2$ and $\theta$ order-preserving. 

Next we show that for any $a \in {\mathbb Q}_{red}$, $R(a)$ is infinite (and similarly for its complement, and for $R(b)$ for $b \in {\mathbb Q}_{blue}$).
If not, let $R(a) = \{b_1, \ldots, b_n\}$ where $b_1 < \ldots < b_n$. Since $|R(a)| \neq 1$, $n \ge 2$. Choose $b_{n+2} > b_{n+1} > b_n$ in 
${\mathbb Q}_{blue}$. By homogeneity there are automorphisms $\theta$, $\varphi$ taking $\{b_1, \ldots, b_n\}$ to $\{b_1, \ldots, b_{n-1}, b_{n+1}\}$
and $\{b_1, \ldots, b_{n-2}, b_{n+1}, b_{n+2}\}$ respectively. Let $\theta a = a_1$ and $\varphi a = a_2$. Since $a$, $a_1$, and $a_2$ have 
different neighbour sets, they must be distinct. 

Suppose first that $a$ is the least or greatest of $a$, $a_1$, $a_2$. Then $p(a) = a$, $p(a_1) = a_2$ is a finite partial automorphism, so
extends to an automorphism $\psi$. As $\psi(a) = a$, $\psi\{b_1, \ldots, b_n\} = \{b_1, \ldots, b_n\}$, and as $\psi$ is order-preserving, 
$\psi b_i = b_i$ for each $i$. But also $\psi \{b_1, \ldots, b_{n-1}, b_{n+1}\} = \{b_1, \ldots, b_{n-2}, b_{n+1}, b_{n+2}\}$, which is impossible 
since $\psi b_{n-1} = b_{n-1}$. Otherwise, $a$ lies between $a_1$ and $a_2$, so we may fix $a_1$ by an automorphism $\psi$ which takes $a$ to $a_2$.
Hence $\psi b_i = b_i$ for $1 \le i \le n-1$ and $i = n+1$. But also $\psi \{b_1, \ldots, b_n\} =  \{b_1, \ldots,b_{n-2}, b_{n+1}, b_{n+2}\}$, so
this is also impossible. This contradiction shows that $R(a)$ is infinite.

Now suppose that for some $a \in {\mathbb Q}_{red}$ there are $b_1 < b_2 < b_3$ in ${\mathbb Q}_{blue}$ such that 
$R(a, b_1) \wedge \neg R(a, b_2) \wedge R(a, b_3)$. Then $p$ given by $p(a) = a$, $p(b_1) = b_3$ is a finite partial automorphism, so extends
to $\theta \in G$. Hence $\{\theta^nb_1: n \in {\mathbb Z}\}$ and  $\{\theta^nb_2: n \in {\mathbb Z}\}$ alternate infinitely between $R(a)$
and ${\mathbb Q}_{blue} \setminus R(a)$. From this it follows that we get a copy of ${\mathbb Q}_2$ in ${\mathbb Q}_{blue}$ by taking the members of $R(a)$
for one `colour', and the members of its complement for the other. A similar argument applies if there are blue $b_1 < b_2 < b_3$ such that
$\neg R(a, b_1) \wedge R(a, b_2) \wedge \neg R(a, b_3)$.

Alternatively, there are no such $b_1, b_2, b_3$. In the first case there are $b_1 < b_2$ such that $R(a, b_1) \wedge \neg R(a, b_2)$. We deduce that
$b \le b_1 \Rightarrow R(a,b)$ and $b \ge b_2 \Rightarrow \neg R(a, b)$, and therefore $R(a) = (-\infty, x)$ or $(-\infty, x]$ for some 
$x \in {\mathbb R}_{blue}$, so ${\mathbb Q}_{blue} \setminus R(a) = [x, \infty)$ or $(x, \infty)$ respectively. 

Next we show that if $R(a) = (-\infty, x)$ then $x$ is irrational. If not, consider $p$ given by $p(a) = a$, $p(x) = x+1$. Then 
$\neg R(a,x) \wedge \neg R(a,x+1)$, so $p$ is a finite partial automorphism. Let $\theta \in G$ extend $p$. 
Then $\theta^{-1}x < x$ so $R(a, \theta^{-1}x)$ but $\neg R(a,x)$, a contradiction.
Next suppose that $R(a) = (-\infty, x]$ where $x$ is rational. This time we let $p(a) = a$, $p(x-1) = x$, and as 
$R(a, x-1)$, $R(a, x)$, $p$ is a partial automorphism. Let $\theta \in G$ extend $p$. Then $\theta x > x$, so $\neg R(a, \theta x)$, again a
contradiction. Similarly if $R(a)$ is of the form $(x, \infty)$ or $[x, \infty)$, $x$ is irrational.

To sum up, there are three possible cases for $R(a)$, namely $(-\infty, x)$, $(x, \infty)$ where $x$ is irrational, or where $R(a)$ and 
${\mathbb Q}_{blue} \setminus R(a)$ form a 2-coloured version of the rationals in ${\mathbb Q}_{blue}$. We need to analyze these a little more,
and also see how the behaviours on ${\mathbb Q}_{red}$ and ${\mathbb Q}_{blue}$ can interact.

First consider the $(x, \infty)$ case for $R(a)$. Let us write $f(a) = x$. Thus $f$ is a function. For any $b \in {\mathbb Q}_{blue}$, 
$R(b)$ and $\neg R(b)$ are non-empty and so there are $c_1, c_2 \in {\mathbb Q}_{red}$ such that $R(c_1,b) \wedge \neg R(c_2, b)$. Since
$G$ acts transitively on the 2-element subsets of ${\mathbb Q}_{red}$, there is $\theta \in G$ taking $\{c_1, c_2\}$ to $\{a_1, a_2\}$.
Then $\theta b$ lies in one of $R(a_1)$, $R(a_2)$ but not the other, and hence $R(a_1) \neq R(a_2)$. Therefore $f(a_1) \neq f(a_2)$. This
shows that $f$ is 1--1.

To see that the image of $f$ is dense in ${\mathbb R}_{blue}$, let $b_1 < b_2$ lie in ${\mathbb Q}_{blue}$. Applying the above argument
the other way round, there is $a \in {\mathbb Q}_{red}$ joined to one of $b_1$, $b_2$ and not the other. We observe that 
$R(a, b_1) \Leftrightarrow b_1 \in R(a) \Leftrightarrow b_1 > f(a)$, and $R(a, b_2) \Leftrightarrow b_2 \in R(a) \Leftrightarrow b_2 > f(a)$.
Since $b_1 < b_2$, we cannot have $b_2 < f(a) < b_1$, and we deduce that $b_1 < f(a) < b_2$. 

To see that the image of $f$ is unbounded above and below, given any $b \in {\mathbb Q}_{blue}$, we note that the image of $f$ intersects
both $(b-1, b)$ and $(b, b+1)$.

Now $R(a)$ is of the form $(f(a), \infty)$ or $(-\infty, f(a))$ for each $a \in {\mathbb Q}_{red}$ (by transitivity of $G$ on ${\mathbb Q}_{red}$
it is the same one of these for each $a$). We can reduce to considering just one of these by applying a `colour reversal'. Given homogeneous
$(X, <, R)$ where $X = {\mathbb Q}_{red} \cup {\mathbb Q}_{blue}$ and ${\mathbb Q}_{red} < {\mathbb Q}_{blue}$, its {\em blue reversal} is got
by reversing the ordering on ${\mathbb Q}_{blue}$, that is, $x \prec y$ if $x,y$ are red and $x < y$, or $x,y$ are blue and $y < x$,
or $x$ is red and $y$ is blue. Then it is easily checked that $(X, \prec, R)$ is homogeneous (essentially using the fact that
${\rm Aut}({\mathbb Q}, <) = {\rm Aut}({\mathbb Q}, >)$). Now for $a \in {\mathbb Q}_{red}$, $R(a)$ is unchanged as a set in passing from
$(X, <, R)$ to $(X, \prec, R)$, but it equals $(f(a), \infty)$ under $<$ if and only if it equals $(-\infty, f(a))$ under $\prec$. This shows
that by applying blue reversal if necessary, we may assume that $R(a) = (f(a), \infty)$ (noting that there is no corresponding change in 
$R(b)$ for $b \in {\mathbb Q}_{blue}$).  

We now show that $f$ is order-preserving or order-reversing. Choose $a_1 < a_2$ in ${\mathbb Q}_{red}$ and suppose that $f(a_1) < f(a_2)$. 
From this we can show that $f$ is order-preserving. (Similarly, if $f(a_2) < f(a_1)$ one shows that $f$ is order-reversing.) Let $a_1' < a_2'$
in ${\mathbb Q}_{red}$ be arbitrary. Then there is an automorphism $\theta$ of $(X, <, R)$ taking $a_1$ to $a_1'$ and $a_2$ to $a_2'$. Since $f$
is definable from $(X, <, R)$, it is preserved by $\theta$, and so
$$f(a_1') = f(\theta a_1) = \theta f(a_1) < \theta f(a_2) = f(\theta a_2) = f(a_2').$$

If we consider the order-reversing case, we apply red reversal of ${\mathbb Q}_{red}$, and the effect is to turn $f$ into an order-preserving map. 

Let us now suppose therefore that $f$ is order-preserving, and for each $a \in {\mathbb Q}_{red}$, $R(a) = (f(a), \infty)$. Let the ordering
$\prec$ be defined on $X = {\mathbb Q}_{red} \cup {\mathbb Q}_{blue}$ by letting $x \prec y$ if $x$ and $y$ have the same colour and $x < y$, 
or $x \in {\mathbb Q}_{red} \wedge y \in {\mathbb Q}_{blue} \wedge f(x) < y$, or $x \in {\mathbb Q}_{blue} \wedge y \in {\mathbb Q}_{red} \wedge x < f(y)$.
Thus the points of ${\mathbb Q}_{red}$ are `slotted in' between those of ${\mathbb Q}_{blue}$, and the fact that $f$ is order-preserving guarantees
that $\prec$ is a linear ordering. To show that the result is isomorphic to ${\mathbb Q}_2$ with the two stated colours it has to be shown that
red and blue points arise in between any two distinct points, and this follows from the facts that both ${\mathbb Q}_{blue}$ and the image of $f$ are 
dense (as subsets of ${\mathbb R}_{blue}$, the order-completion of ${\mathbb Q}_{blue}$). Thus $X$ is bounded generic, as we have essentially reversed
the construction given in the preamble. This gives this case, together with the three other structures obtained by performing colour reversals.

The same argument applies from the assumption that $R(b)$ is bounded (above or below) for each blue point $b$, which therefore gives rise to 
the same list of structures. 

So it remains to treat the case in which each $R(a)$ is dense and codense in ${\mathbb Q}_{blue}$, and each 
$R(b)$ is dense and codense in ${\mathbb Q}_{red}$. We shall deduce from this situation the apparently stronger statement that for each
pair of finite disjoint subsets $A_1$ and $A_2$ of ${\mathbb Q}_{red}$ the set of all points of ${\mathbb Q}_{blue}$ which are adjacent to
all members of $A_1$ and to no members of $A_2$ is dense in ${\mathbb Q}_{blue}$. 

We start by noting that this set is non-empty. Let us enumerate $A_1 \cup A_2$ in increasing order as $a_1 < a_2 < \ldots < a_n$, and we choose 
$a_i'$ inductively in such a way that $a_i \in A_1 \Leftrightarrow a_i' \in R(b)$ and $a_i \in A_2 \Leftrightarrow a_i' \not \in R(b)$, where 
$b$ is some fixed member of ${\mathbb Q}_{blue}$. Since $R(b)$ and ${\mathbb Q}_{red} \setminus R(b)$ are non-empty, we can find 
$a_1' \in R(b)$ if $a_1 \in A_1$ and $a_1' \not \in R(b)$ if $a_1 \in A_2$. A similar argument applies to the induction step, using the
facts that $R(b)$ and ${\mathbb Q}_{red} \setminus R(b)$ are unbounded above. By homogeneity of $(X, <, R)$ there is an automorphism $\theta$
of $(X, <, R)$ which takes $a_i'$ to $a_i$ for each $i$. Let $b' = \theta b$. Then if $a_i \in A_1$, $a_i' \in R(b)$, so 
$a_i = \theta a_i' \in R(\theta b) = R(b')$, and if $a_i \in A_2$, $a_i' \not \in R(b)$, so $a_i \not \in R(b')$. Hence $b'$ lies in
the desired set.

We can now show by induction on $|A|$ that for any disjoint $A_1$ and $A_2$ such that $A = A_1 \cup A_2$, the set of possible witnesses 
corresponding to that choice of $A_1$ and $A_2$ (values of $b'$) is dense in ${\mathbb Q}_{blue}$. The case $|A| = 1$ is the assumption. 
Otherwise consider a larger value $n$ of $|A|$, and suppose inductively that the result holds for smaller values. In fact we treat all
possible subdivisions of $A$ into sets $A_1$ and $A_2$ simultaneously. To handle this, we let $P$ be the set of all maps from
$A = \{a_1, \ldots, a_n\}$ (enumerated in increasing order as $\bf a$) into $\{+, -\}$, and the partition corresponding to $p \in P$ is then
$A_1 = \{i: p(a_i) = +\}$ and $A_2 = \{i: p(a_i) = -\}$. For each $p \in P$ we let 
$W_p = \{b \in {\mathbb Q}_{blue}: (\forall i)(p(a_i) = + \Rightarrow R(a_i, b) \wedge p(a_i) = - \Rightarrow \neg R(a_i, b))\}$. Members 
of $W_p$ are thought of as `witnesses' corresponding to $p$ (and $\bf a$). Then the previous paragraph showed that each $W_p$ is non-empty (and
we are trying to show that it is dense). 

We now form a new `colouring' of the points of ${\mathbb Q}_{blue}$ in which two points receive the same colour provided that the sets of 
members of $A$ that they are joined to are equal. This is equivalent to saying that they are witnesses for the same $p$. Under this 
colouring, we can see that ${\mathbb Q}_{blue}$ becomes a homogeneous coloured linearly ordered set. For if $q$ is a finite partial automorphism, 
we can extend $q$ to act on $A \cup {\rm dom} \, q$ by fixing each member of $A$, and this extension clearly preserves the bipartite relation, 
so is a finite partial automorphism of $(X, <, R)$. By homogeneity this extends to $\theta \in G$, and the restriction of $\theta$ to 
${\mathbb Q}_{blue}$ extends $q$, and also preserves the new `colours'. According to the classification of (finitely coloured) countable 
homogeneous coloured linear orders (see for instance \cite{Campero}), ${\mathbb Q}_{blue}$ is the finite union of convex pieces, on 
each of which the `colours' which arise occur densely. By induction hypothesis, the witnesses corresponding to any $p$ on a set of size 
$n-1$ occur densely. These are precisely the same as the union of the sets of witnesses of each of the two possible extensions of $p$ to 
a domain of size $n$. We deduce that ${\mathbb Q}_{blue}$ is cut into at most two pieces. If there is only one piece, then the set of
witnesses $W_p$ corresponding to this value of $p$ for every $p \in P$ is dense, and this gives what is required. We note that the notation 
$W_p$ tacitly includes reference to $\bf a$, since this is the domain of $p$.

So now we suppose for a contradiction that there are (exactly) two pieces, which must have the forms $(-\infty, x)$ and $(x, \infty)$ for some 
irrational $x$. The argument just given shows that the witnesses for any two values of $p$ which differ on just one place are in different
pieces, and one easily deduces (from the fact that the graph structure on $P$ given by saying that two such maps are joined provided that they 
differ in exactly one position is connected) that one of the pieces contains all the witnesses for those $p$ which have an even number of $+$ 
values, and the other contains all the witnesses for those $p$ which have an odd number of $+$ values.

Following the method for the basis case, we can define $f(a_1, \ldots, a_n) = x$. We now see that $f$ is 1--1 on each co-ordinate. That is, if 
$\bf a$ and $\bf a'$ are sequences differing in exactly one position, the $i$th say, then $f({\bf a}) \neq f({\bf a'})$. For this, take
any $p \in P$, and extend it to $q$ defined on $A \cup \{a_i'\}$ by letting $q(a_i') \neq p(a_i)$. Then $q$ is a map on $n+1$ points,
still to $\{+, -\}$, and so we know that there is a witness $b$ for $q$ in ${\mathbb Q}_{blue}$. Since $q$ extends $p$, $b \in W_p$, and
if we let $p'$ be the restriction of $q$ to ${\bf a}'$, then similarly, $b \in W_{p'}$. If $f({\bf a}) = f({\bf a}') = x$, then $W_p$ and $W_{p'}$
are complementary (one equal to $(-\infty, x)$, the other to $(x, \infty)$), so they cannot both contain $b$, giving a contradiction.

Extending this we can see that $f$ is either order-preserving or order-reversing on each co-ordinate. For as just shown if 
$\bf a$ and $\bf a'$ are sequences differing in the $i$th position, then $f({\bf a}) \neq f({\bf a'})$, so either
$f({\bf a}) < f({\bf a'})$ or $f({\bf a}) > f({\bf a'})$, suppose the former for instance. Then by homogeneity, we can take 
$\{a_j: 1 \le j \le n\} \cup \{a_i'\}$ to any other similar $(n+1)$-element subset of ${\mathbb Q}_{red}$, and so the ordering must be the
same way round (either $<$ or $>$).

First treat the case in which $f({\bf a}) < f({\bf a'})$, and choose $b < f({\bf a})$. Consider the relations between $a_j$ and $b$. Since
${\bf a}'$ differs from $\bf a$ in just the $i$th place, the relations between $a_j'$ and $b$ are the same except possibly between $a_i'$
and $b$. However, $b < f({\bf a}')$ too, so there cannot be a difference on just one place, so $R(a_i, b) \Leftrightarrow R(a_i', b)$.
Since $R(b)$ and its complement are dense in ${\mathbb Q}_{red}$, there must be some $a_i' \in (a_i, a_{i+1})$ such that 
$R(a_i, b) \Leftrightarrow \neg R(a_i', b)$, which gives a contradiction.

If instead $f({\bf a}) > f({\bf a'})$, we choose $b > f({\bf a})$, and we consider the relations between $a_j$ and $b$. Again, since ${\bf a}'$ 
differs from $\bf a$ in just the $i$th place, the relations between $a_j'$ and $b$ are the same except possibly between $a_i'$ and $b$. This time, 
$b > f({\bf a}')$, so as before, $R(a_i, b) \Leftrightarrow R(a_i', b)$ and we reach a contradiction as before.

The conclusion is that each $W_p$ is dense. From this it follows that $(X, <, R)$ is unbounded generic. For if $A_1$ and $A_2$ are finite
disjoint subsets of ${\mathbb Q}_{red}$, we can let $p \in P$  be given by $p(a_i)= +$ if $a_i \in A_1$ and $p(a_i) = -$ if $a_i \in A_2$.
Since as we have shown, $W_p$ is dense in ${\mathbb Q}_{blue}$, it follows that so is the set of all members of ${\mathbb Q}_{blue}$
that are $R$-related to all members of $A_1$ and to none of $A_2$. A standard back-and-forth argument now shows that $(X, <, R)$ is unbounded 
generic (and is unique up to isomorphism).

(v)  We consider $({\mathbb Q}.2, <)$ in which without loss of generality, each pair is coloured (red, blue). For a red point $a$ paired with the blue
point $b$, we see what $R(a)$ can be. If $b_2 > b_1 > b$ are all blue, then $p$ given by $p(b) = b$, $p(b_1) = b_2$ is a finite 
partial automorphism, so extends to $\theta \in G$. Since $\theta(b) = b$, also $\theta(a) = a$. Hence $R(a, b_1) \Leftrightarrow R(a, b_2)$.
This shows that $b_1 \in R(a) \Leftrightarrow b_2 \in R(a)$. Similarly if $b_1 < b_2 < b$. We deduce that $R(a) = \{b\}$ or 
$(b, \infty) \cap {\mathbb Q}_{blue}$ or $(-\infty, b) \cap {\mathbb Q}_{blue}$, or a union of these. Clearly whichever option holds 
will be the same for all $a \in {\mathbb Q}_{red}$. So
this gives 8 possibilities, as stated in the theorem. One checks easily that all of these are homogeneous.

(vi) The underlying ordered set is ${\mathbb Q}_2$ in this case. It is easy to check that all the listed bipartite graphs are homogeneous.
Let us see that they are the only possibilities. We follow some of the same steps as in case (iv). 

We assume that $R$ is not empty or complete (which are vacuously unions of the stated structures). We concentrate on analyzing 
$R_1 = \{(a, b) \in R: a \mbox{ red}, b \mbox{ blue}, a < b\}$, and suppose that $R_1$ is not empty or right complete. Hence for some red $a$, 
$R_1(a)$ does not equal $(a, \infty) \cap {\mathbb Q}_{blue}$. In the first place we note that $R_1(a)$ cannot be finite non-empty. For 
if so, let it be $\{b_1, \ldots, b_n\}$, and choose blue $b' \in (b_{n-1}, b_n)$ and $b'' > b_n$. Then the map which fixes $a$ and sends 
$b'$ to $b''$ is a finite partial automorphism, but it cannot extend to an automorphism $\theta$ because $\theta$ must fix 
$\{b_1, \ldots, b_n\}$ pointwise, contrary to being order-preserving. Similarly, 
the complement $\neg R_1(a)$ of $R_1(a)$ in $(a, \infty) \cap {\mathbb Q}_{blue}$ is also infinite.

Next we see that each of $R_1(a)$ and $\neg R_1(a)$ is dense in $(a, \infty)$. If $b_1 < b_2 < b_3$ are blue and greater than $a$, and 
$b_1, b_3 \in R_1(a)$ and $b_2 \not \in R_1(a)$, then there is an automorphism fixing $a$ and taking $b_1$ to $b_3$, so there is an 
infinite alternation between $R_1(a)$ and its complement in $(a, \infty)$. Similarly if $b_1, b_3 \not \in R_1(a)$ and $b_2 \in R_1(a)$. 
Hence in each of these cases there are blue $b_1 < b_2 < b_3$ with $b_1, b_3 \in R_1(a)$ and $b_2 \not \in R_1(a)$. Then we can take 
any blue members $b_1' < b_3'$ of $R_1(a)$ to $b_1$ and $b_3$ by an automorphism fixing $a$, so there is a member of $\neg R_1(a)$ 
between $b_1'$ and $b_3'$. The same proof as in case (iv) shows that each of $R_1(a)$ and $\neg R_1(a)$ is dense in $(a, \infty)$. 

If on the other hand there are no such $b_1, b_2, b_3$ we deduce that $R_1(a)$ is either an initial or final segment of the set of
blue members of $(a, \infty)$. Let $x$ be the cut between $R_1(a)$ and its complement. By homogeneity, there is an automorphism that fixes 
$a$ and maps a red member of $(a, x)$ to a red member of $(x, \infty)$. However this must fix $x$, contrary to its being order-preserving.
We deduce that each of $R_1(a)$ and $\neg R_1(a)$ is dense in $(a, \infty)$. 

Next we show that for any finite disjoint sets $A_1$ and $A_2$  of red points, the set of blue points which are joined to all members of 
$A_1$ and to none of $A_2$ is dense in $(\max (A_1 \cup A_2), \infty)$. We follow essentially the same method as for (iv). The proof
goes by induction on $|A_1 \cup A_2|$, the case where this is 1 having been already covered. Write $A = A_1 \cup A_2$, and $P$ for the
family of maps from $A$ to $\{+, -\}$. The case for $A_1$ and $A_2$ is obtained by letting $p(a) = +$ if $a \in A_1$ and $p(a) = -$
if $a \in A_2$, and we have the notion of a `witness' corresponding to $p$ as before, with $W_p$ being the set of all witnesses.

Following the proof from (iv), first note that each $W_p$ is non-empty. Next we use induction on $|A|$ to show that each $W_p$ is 
dense in $(\max (A), \infty)$. Form a new colouring of the blue points of $(\max (A), \infty)$ in which two points receive the same 
colour provided that the sets of members of $A$ that they are joined to are equal. Then $(\max (A), \infty) \cap {\mathbb Q}_{blue}$ becomes
a homogeneous coloured linearly ordered set. Thus $(\max (A), \infty)$ is the finite union of convex pieces, on each of which the 
`colours' which arise occur densely. By induction hypothesis, the witnesses corresponding to any $p$ on a set of size 
$n-1$ occur densely. These are precisely the same as the union of the sets of witnesses of each of the two possible extensions of $p$ to 
a domain of size $n$. We deduce that $(\max (A), \infty)$ is cut into at most two pieces. If there is only one piece, then the set of
witnesses $W_p$ corresponding to this value of $p$ for every $p \in P$ is dense, and this gives what is required. 

So now we suppose for a contradiction that there are (exactly) two pieces, which must have the forms $(\max (A), x)$ and $(x, \infty)$ for 
some irrational $x$. The argument just given shows that the witnesses for any two values of $p$ which differ on just one place are in different
pieces, and hence one of the pieces contains all the witnesses for those $p$ which have an even number of $+$ values, and the other contains 
all the witnesses for those $p$ which have an odd number of $+$ values.

Now write $x$ as $f(a_1, \ldots, a_n)$ where this is the increasing enumeration of $A$, and as before one shows that $f$ is 1--1 on each
co-ordinate, and indeed is either order-preserving or order-reversing on each co-ordinate. Now consider sequences
$\bf a$ and $\bf a'$ which differ in the $i$th position, and first treat the case in which $f({\bf a}) < f({\bf a'})$, and choose 
$b \in (\max (A), f({\bf a}))$. It follows that also $b \in (\max (A), f({\bf a}'))$. Since ${\bf a}'$ differs from $\bf a$ in just 
the $i$th place, the relations between $a_j'$ and $b$ are the same except possibly between $a_i'$ and $b$, and as there cannot be 
a difference on just one place, $R_1(a_i, b) \Leftrightarrow R_1(a_i', b)$. Since however $R_1(b)$ and its complement are dense in 
$(-\infty, b)$, there must be some $a_i' \in (a_i, a_{i+1})$ such that $R_1(a_i, b) \Leftrightarrow \neg R_1(a_i', b)$, which gives a contradiction.

If instead $f({\bf a}) > f({\bf a'})$, we choose $b > f({\bf a})$, and apply a similar argument.

The conclusion is that each $W_p$ is dense. From this it follows that $(X, <, R_1)$ is unbounded generic. For if $A_1$ and $A_2$ are finite
disjoint subsets of ${\mathbb Q}_{red}$, we can let $p \in P$  be given by $p(a_i)= +$ if $a_i \in A_1$ and $p(a_i) = -$ if $a_i \in A_2$.
Since as we have shown, each $W_p$ is dense in $(\max (A), \infty) \cap {\mathbb Q}_{blue}$, it follows that so is the set 
of all members of $(\max (A), \infty) \cap {\mathbb Q}_{blue}$ that are $R_1$-related to all members of $A_1$ and to none of $A_2$. A standard 
back-and-forth argument now shows that $(X, <, R_1)$ is right generic.

We have shown that $R_1$ is empty, right generic, or right complete. Similarly $R_2 = \{(a, b) \in R: a \mbox{ red}, b \mbox{ blue}, b < a\}$
is empty, left generic, or left complete, and this completes the proof.   $\Box$

\vspace{.2in}

\noindent{\bf Concluding remarks}

If we wanted to extend the results presented to multipartite graphs, it would be necessary to combine the techniques given here with those
of \cite{Jenkinson}. The bipartite restrictions should all be homogeneous, hence in the list given here, but one would have to analyze
how they can fit together. In \cite{Jenkinson} this aspect was quite straightforward, as there were so few possibilities for bipartite restrictions 
(just one non-trivial one, generic), and one could focus on which multipartite graphs could be omitted. It was shown that any minimally
omitted graph was `monic' (having at most one vertex in each part). A specially interesting case was the existence of so-called
`omission quartets', being a configuration in the quadripartite case, of four monic graphs each on three vertices, which could be
jointly omitted. The main burden of the remainder of the proof was that on any greater number of parts, the situation is still 
controlled by individual monics, and omission quartets. In the ordered case, the complications promise to be considerably more
involved.

\end{document}